\documentclass[12pt]{article}

\usepackage[english]{babel}
\usepackage{amssymb}
\usepackage{graphics}

\newcommand{\qed}{$\;\;\;\Box$}

\newcounter{claim}[section]
\renewcommand{\theclaim}{\arabic{claim}}
{\par\medskip\par}

{\qed\par\smallbreak}

\newcommand{\hide}[1]{}
\newtheorem{thm}{Theorem}[section]

\newtheorem{defn}[thm]{Definition}

\date{}
\title{Smoothing property of solutions to \\ nonlocal hyperbolic problems}

\newcounter{thesame}
\author{
	Iryna Kmit
	\thanks{Institute of Mathematics, Humboldt University of Berlin. On leave from the
		Institute for Applied Problems of Mechanics and Mathematics, 
		National Academy of Sciences of	Ukraine, Lviv, Ukraine. 
		{\small   E-mail:
			{\tt kmit@mathematik.hu-berlin.de}}}
}

\begin{document}
\maketitle

\begin{abstract}
	\noindent
We consider nonlocal initial boundary value problems with  integral boundary conditions for integro-differential 
first order hyperbolic systems. We prove a 
general regularity result stating that the $L^2$-generalized solutions become eventually  continuous.

\end{abstract}

\section{Problem setting  and motivation}
\renewcommand{\theequation}{{\thesection}.\arabic{equation}}
\setcounter{equation}{0}

\subsection{Statement of the problem and our result}

We will consider  integro-differential 
first order hyperbolic system of the type 
\begin{equation}
 \begin{array}{ll}
 \displaystyle\partial_{t}u_j
 +a_j(x,t)\partial_{x}u_j
 +\sum_{k=1}^{n}b_{jk}(x,t)u_k
 +\sum_{k=1}^{n}\int_{0}^{x}g_{jk}(y,t)u_k(y,t)dy= f_j(x,t),\\
\qquad\qquad\qquad\qquad\qquad\qquad\qquad\qquad\qquad\;\;\; (x,t)\in(0,1)\times(0,\infty),\; j\le n,&
 \end{array}
 \label{eq:1}
 \end{equation}
 subjected to the initial conditions  
 \begin{equation}\label{eq:2}
 u_{j}(x,0)=\varphi_{j}(x), \quad x\in[0,1],\;\;\; j\le n, 
 \end{equation}
 and the integral boundary conditions 
 \begin{equation}\label{eq:3}
 \begin{array}{ll}
   u_{j}(0,t)=
 \displaystyle\sum_{k=1}^{n}\int_{0}^{1}r_{jk}(x,t)u_{k} dx,\quad t\in[0,\infty),
 \;\;\; 1\le j\le m,  \\ [2mm]
  u_{j}(1,t)=
 \displaystyle\sum_{k=1}^{n}\int_{0}^{1}r_{jk}(x,t)u_{k} dx,,\quad t\in[0,\infty),
 \;\;\; m<j\le n,
 \end{array}
 \end{equation}
where $0\le m\le n$ are fixed integers and $n\ge 2$.
In this form, which is  motivated by 
applications,  the problem has been studied in \cite{KrsSm,SanNak}.
First order hyperbolic systems with smoothing boundary conditions of the integral type
 (\ref{eq:3})
 appear, in particular, in  applications to  
  population dynamics \cite{eft,mogulru,webb}.

It is known that, whatsoever $\varphi_i\in C\left([0,1]\right)$, 
there exists a unique peacewise
continuous solution to the initial-boundary value problem (\ref{eq:1})--(\ref{eq:3})
 with possible first order discontinuities along
characteristic curves emanating from the points $(0,0)$ and $(1,0)$. 
The discontinuities dissappear only whenever the zero order compatibility condition  between 
(\ref{eq:3}) and (\ref{eq:2})
  is fulfilled. It turns out that for some classes of boundary conditions one can speak about 
the disappearance of the discontinuities {\it after} some time even if the zero order compatibility 
conditions at points $(0,\tau)$ and $(1,\tau)$ are not fulfilled. 
Even more, as it follows from \cite{Km,KL}, one can expect that
for linear problems this kind of smoothing property (higher regularity of solutions after
a ``smoothing time'') does not depend on the regularity of the initial data.
Our aim is to show that generalized solutions to the problem  (\ref{eq:1})--(\ref{eq:3})
become eventually continuous, despite of the initial data are supposed to be only $L^2$-regular.

Set 
$$
\Omega=\{(x,t)\,:\,0< x< 1, 0<t<\infty\}.
$$
Suppose that 
\begin{equation}\label{eq:L1}
\inf_{(x,t)\in \overline\Omega}a_j>0 \mbox{ for all } j\le m\quad \mbox{ and }\quad 
\sup_{(x,t)\in \overline\Omega}a_j<0 \mbox{ for all } j>m.
\end{equation}
This entails   that the system (\ref{eq:1}) is non-degenerate.
Let us make a couple of
  regularity assumptions on the coefficients of (\ref{eq:1}), namely
\begin{equation}\label{eq:L3}
a_j,  b_{jk}  \in C^1(\overline\Omega)  \mbox{ and  } g_{jk}, r_{jk}\in C^{0,1}_{x,t}(\overline\Omega) 
\mbox{ for all }
j,k\le n
\end{equation}
and
\begin{equation}\label{cass1}
\begin{array}{l}
\mbox{for all } 1 \le j \not= k \le n \mbox{ there exists }  \beta_{jk} \in C^1
\\\mbox{such that } 
b_{jk}=\beta_{jk}(a_k-a_j).
\end{array}
\end{equation}
It turns out that the last property is crucial for our analysis.

By $L^2(0,1)^n$ we  denote the vector space of  functions  $u=(u_1,\dots,u_n)$
with  $u_j\in L^2(0,1)$. 

Let us introduce the notion of an $L^2$-generalized solution to the problem (\ref{eq:1})--(\ref{eq:3}). 
Fix  arbitrary  $\varphi_j\in L^2(0,1)$ and 
sequences  $\varphi_j^l\in C_0^\infty([0,1])$ 
such that   $\varphi_j^l\to\varphi_j$ in $L^2(0,1)$. Note that, due to the fact that  $\varphi_j^l$ are compactly supported for all 
$l\in\mathbb{N}$ and $j\le n$, they satisfy 
the zero order and the first order compatibility conditions  between (\ref{eq:2}) and  (\ref{eq:3}). By \cite{ijdsde,lakra},
given  $l\in\mathbb{N}$, the problem (\ref{eq:1})--(\ref{eq:3}) has a unique classical solution, say $u^l$.
Similarly to \cite{KL} one can prove that  these solutions satisfy the following estimate:
\begin{equation}\label{bound40}
\max_{j\le n}\|u^l(\cdot,t)\|_{L^2(0,1)}\le 
M e^{\omega t}\max_{j\le n}\|\varphi_j\|_{L^2(0,1)}\quad \mbox{ for all } t>0,
\end{equation}  
for some constants $M,\omega$ not depending on $t$, $\varphi_j$,   and $l\in\mathbb{N}$. Hence, there exist
a unique vector-function  $ u \in C([0,\infty), L^2(0,1))^n$ such that 
$$
 \|u(\cdot,\theta)-u^{l}(\cdot,\theta)\|_{L^2(0,1)^n}\to 0\quad 
\mbox{ as } l\to\infty,
$$
uniformly in $\theta$ varying in the range $0\le\theta\le t$, for every $t>0$. 
The vector-function $u$ is called an  $L^2$- {generalized  solution} 
to the problem (\ref{eq:1})--(\ref{eq:3}).
Furthermore, the following estimate is true:
\begin{equation}\label{bound4}
\|u(\cdot,t)\|_{L^2(0,1)^n}\le 
M e^{\omega t}\max_{j\le n}\|\varphi_j\|_{L^2(0,1)}\quad  \mbox{ for all } t>0,
\end{equation} 
what follows from (\ref{bound40}).

Our main result states that all $L^2$-generalized solutions have a smoothing property, in the following sense.

\begin{defn}\label{defn:smoothing_abstr}  
The $L^2$-generalized solution to the problem (\ref{eq:1})--(\ref{eq:3})
is called {\it eventually smoothing} if there is $d>0$ such that 
$u$ is a continuous vector-function after the time  $t\ge d$.
\end{defn}

This property means that the solutions become more regular {\it in a finite time}
if to compare with their
regularity {\it in the entire domain}. This  contrasts to the parabolic case where, due to the infinite propagation speed, the smoothing 
property for solutions is encountered  in the entire domain (immediately after leaving the 
initial axis).

\begin{thm}\label{main_smoothing}
Suppose that the conditions   (\ref{eq:L1}),  
( \ref{eq:L3}) and (\ref{cass1})  are fulfilled. Then the 
$L^2$-generalized solution to the problem (\ref{eq:1})--(\ref{eq:3})
is eventually smoothing in the sense of Definition \ref{defn:smoothing_abstr}.
\end{thm}

\subsection{Integral representation of the problem (\ref{eq:1})--(\ref{eq:3})}

 For given $j\le n$,  $x \in [0,1]$, and $t \in \mathbb{R}$, the $j$-th characteristic of (\ref{eq:1})
passing through the point $(x,t)$ is defined
as the solution $\xi\in [0,1] \mapsto \omega_j(\xi,x,t)\in \mathbb{R}$ to the initial value problem
\begin{equation}\label{char}
\partial_\xi\omega_j(\xi,x,t)=\frac{1}{a_j(\xi,\omega_j(\xi,x,t))},\;\;
\omega_j(x,x,t)=t.
\end{equation} 
Due to (\ref{eq:L1}), the characteristic curve $\theta=\omega_j(\xi,x,t)$ reaches the
boundary of $\Omega$ in two points with distinct ordinates. Let $x_j=x_j(x,t)$
denote the abscissa of that point whose ordinate is smaller.
Note that  (\ref{eq:L1}) entails that $x_j$ does not depend on $x,t$ but only on $j$
and, therefore, is given by the formula
$$
x_j=\left\{
\begin{array}{rcl}
0 & \mbox{for} &  1\le j\le m,\\
1 & \mbox{for} & m< j\le n.
\end{array}
\right.
$$
Further we will simply write $\omega_j(\xi)$ for $\omega_j(\xi,x,t)$.
Write
\begin{equation}\label{cd}
\begin{array}{cc}
\displaystyle c_j(\xi,x,t)=\exp \int_x^\xi
\left(\frac{b_{jj}}{a_{j}}\right)(\eta,\omega_j(\eta))\,d\eta,\quad
d_j(\xi,x,t)=\frac{c_j(\xi,x,t)}{a_j(\xi,\omega_j(\xi))}.
\end{array}
\end{equation} 
Let us introduce   linear bounded operators $R,B,G: C(\overline\Omega)^n \mapsto C(\overline\Omega)^n$  by
\begin{eqnarray}
 (Ru)_j(x,t)&=&
 c_j(x_j,x,t)\sum_{k=1}^{n}\int_{0}^{1}r_{jk}(\eta,\omega_j(x_j))u_{k}(\eta,\omega_{j}(x_j))\,d\eta, \nonumber
\\
 (Bu)_j(x,t)&=&
 -\sum_{k\neq j}\int_{x_j}^{x}d_j(\xi,x,t)b_{jk}(\xi,\omega_j(\xi))u_k(\xi,\omega_j(\xi))\, d\xi,
 \label{f34} \\
(Gu)_j(x,t)&=& 
-\sum_{k=1}^{n}\int_{x_j}^{x}\int_{0}^{\xi}d_{j}(\xi,x,t)g_{jk}(y,\omega_j(\xi))
u_{k}(y,\omega_j(\xi))\, dyd\xi, \label{f3014}
 \end{eqnarray}

Straightforward calculations  show that a $C^1$-map $u:\overline{\Omega} \to \mathbb{R}^n$ 
is a solution to 
the PDE problem (\ref{eq:1})--(\ref{eq:3}) if and only if
it satisfies the following system of integral equations
\begin{eqnarray}
\label{rep1}
u_j(x,t)=
\left(Su\right)_j(x,t)
-\int_{x_j}^x d_j(\xi,x,t)\sum_{k\not=j} (b_{jk}u_k)(\xi,\omega_j(\xi))\,d\xi,
\quad j\le n,
\end{eqnarray}
where the affine bounded operator $S$
is defined by
\begin{eqnarray}
\label{Q}
(Su)_j(x,t)=
\left\{\begin{array}{lcl}
\displaystyle  \left(Ru\right)_j(x,t)
& \mbox{if}&  x_j=0 \mbox{ or } x_j=1 \\
c_j(x_j,x,t)\varphi_j(x_j)      & \mbox{if}&  x_j\in(0,1)
\end{array}\right.
\end{eqnarray}
on a subspace of $C(\overline\Omega)^n$ of functions satisfying  (\ref{eq:2}).

\section{Proof of Theorem \ref{main_smoothing}}\label{sec:smooth}

Here we prove Theorem \ref{main_smoothing}. To this end, we will use 
results about existence and uniqueness of continuous and classical solutions to the problem (\ref{eq:1})--(\ref{eq:3})
 proved in \cite{ijdsde,Km1}.

The proof extends the ideas of
\cite{kmit,Km}
where
 the smoothing property is proved for the initial data in 
the space of continuous functions satisfying the zero order compatibility conditions. 
In \cite{kmit,Km} we  
show that the solutions reach the $C^k$-regularity in a finite time for each $k$.
Here we follow a similar argument and extend the smoothing results to the case where the initial data 
are $L^2$-functions only.

Our starting point is that, given $\varphi\in L^2(0,1)^n$, the initial boundary value
 problem (\ref{eq:1})--(\ref{eq:3}) has a unique solution 
$u\in C\left(\mathbb{R},L^2(0,1)\right)^n$. Our aim is to show that 
$u(x,t)\in C\left(\overline\Omega_{d}\right)^n$ for some $d>0$, where
$$
\Omega_\tau=\{(x,t)\,:\,0< x< 1, \tau<t<\infty\},\quad\tau>0.
$$
Let $d$ be chosen so that 
$d\ge\omega_j(1,\omega_k(1,0,0),0)$ for all $j,k\le m$ and $d\ge\omega_j(0,\omega_k(0,1,0),1)$ for all $m+1\le j,k\le n$.
Fix an arbitrary $\varphi\in L^2(0,1)^n$ and consider 
 an arbitrary sequence $\varphi^l\in C([0,1])^n$  of functions satisfying the zero order
compatibility conditions between (\ref{eq:3}) and~(\ref{eq:2}) 
such that $\varphi^l\to\varphi$ in $L^2(0,1)^n$ (one can easily prove that the sequence $\varphi^l$ exists). 
Denote by $u^l$ the continuous solution to
(\ref{eq:1})--(\ref{eq:3}) with $\varphi$ replaced by  $\varphi^l$.
Moreover, we have (see \cite[Lemma 4.2]{KL})
\begin{equation}\label{seq1}
u^l\to u \mbox{ in }  C\left([\tau,\theta],L^2(0,1)\right)^n \mbox{ as } l\to\infty,
\end{equation}
for any $\theta>0$, where $u\in C\left([0,\theta],L^2(0,1)\right)^n$ is the generalized 
$L^2$-solution to
the problem (\ref{eq:1})--(\ref{eq:3}). To prove the theorem, it is sufficient to
 show that 
\begin{equation}\label{seq2}
u^l \mbox{ converges in }  C\left(\overline\Omega_{4d}^{4d+\alpha}\right)^n
\mbox{ as } l\to\infty,
\end{equation}
for any $\alpha>0$. Here and below, given $\beta<\gamma$, we write
$$
\Omega_\beta^\gamma=\Omega_{\beta}\setminus\overline\Omega_{\gamma}.
$$

For the solution 
$u^l\in  C\left(\overline\Omega\right)^n$ 
restricted to $\overline\Omega_{d}$ 
we have the operator representation:
\begin{eqnarray}
\label{rep_u0}
u^l\big|_{\overline\Omega_{d}}=Bu^l+Gu^l+Ru^l.
\end{eqnarray}
On the account of the choice of $d$, the function $u^l$ in the right-hand side of (\ref{rep_u0}) fulfills the same equation 
(\ref{rep_u0}),
what  entails
\begin{eqnarray*}
u^l\big|_{\overline\Omega_{d}}=(B+R)(B+R+G)u^l+Gu^l.
\end{eqnarray*}
We are done if we prove that the right-hand side of the last equality converges in 
$C\left(\overline\Omega_{2d}^{2d+\alpha}\right)^n$ for any $\alpha>0$
as $l\to\infty$.
Fix an arbitrary $\alpha>0$.
Since the operators $B$ and $R$ are bounded, it suffices to prove that the 
sequences
\begin{equation}\label{conv_all}
 B^2u^l, RBu^l, BRu^l, R^2u^l, Gu^l
\mbox{ converge in } C\left(\overline\Omega_{2d}^{2d+\alpha}\right)^n
\end{equation}
 as  $l\to\infty$.

We start with $\left(B^2u^l\right)(x,t)$.
Given $j\le n$, consider the following expression for $\left(B^2u^l\right)_j(x,t)$,
obtained after changing the order of integration:
\begin{equation}\label{D11}
\begin{array}{cc}
\left(B^2u^l\right)_j(x,t)
=\displaystyle\sum_{k\not=j}\sum_{i\not=k}
\int_{x_j}^x \int_\eta^x d_{jki}(\xi,\eta,x,t)b_{jk}(\xi,\omega_j(\xi))\\
\displaystyle\times
u_i^l(\eta,\omega_k(\eta,\xi,\omega_j(\xi))) d \xi d \eta
\nonumber
\end{array}
\end{equation}
with
\begin{eqnarray*}
d_{jki}(\xi,\eta,x,t)
=d_j(\xi,x,t)d_k(\eta,\xi,\omega_j(\xi))b_{ki}(\eta,\omega_k(\eta,\xi,\omega_j(\xi))).
\label{djkl}
\end{eqnarray*}
Fix $k\not=j$.
Let us change the variables 
\begin{equation}\label{change}
\xi\in[0,1]\mapsto\theta=\omega_k(\eta,\xi,\omega_j(\xi)).
\end{equation}
Taking into the account 
the formulas
\begin{eqnarray}
\label{dx}
\partial_x\omega_j(\xi,x,t) & = & -\frac{1}{a_j(x,t)} \exp \int_\xi^x 
\frac{\partial_ta_j(\eta,\omega_j(\eta))}{a_j(\eta,\omega_j(\eta))^2}\, d\eta,\\
\label{dt}
\partial_t\omega_j(\xi,x,t) & = & \exp \int_\xi^x 
\frac{\partial_ta_j(\eta,\omega_j(\eta))}{a_j(\eta,\omega_j(\eta))^2}\, d\eta,
\end{eqnarray}
from (\ref{change}) we get
\begin{eqnarray}\label{theta}
d\theta& = & \left[\partial_2\omega_k(\eta,\xi,\omega_j(\xi))+\partial_3\omega_k(\eta,\xi,\omega_j(\xi))\partial_\xi\omega_j(\xi)\right]
d\xi\nonumber\\
\displaystyle
& = &\frac{a_k(\xi,\omega_j(\xi))-a_j(\xi,\omega_j(\xi))}{a_j(\xi,\omega_j(\xi))a_k(\xi,\omega_j(\xi))}
\partial_3\omega_k(\eta,\xi,\omega_j(\xi)) d\xi,
\end{eqnarray}
where $\partial_k$ denotes the partial derivative with respect to the $k$-th argument.
It follows that (\ref{change}) is non-degenerate  for all $\xi\in[0,1]$
fulfilling the condition   $a_k(\xi,\omega_j(\xi))-a_j(\xi,\omega_j(\xi))\ne 0$.
The inverse of (\ref{change}) for those $\xi$ 
will be denoted by $\tilde x(\theta,\eta,x,t)$. Moreover, we have the following identity:
\begin{equation}\label{kj}
\omega_k(\tilde x(\theta,\eta,x,t),\eta,\theta)=\omega_j(\tilde x(\theta,\eta,x,t),x,t).\nonumber
\end{equation}
Therefore, after changing the variables (\ref{change})
we  come up with the following formula for the summands contributing into $\left(B^2u^l\right)$:
\begin{equation}\label{after_change}
\begin{array}{cc}
\displaystyle\int_{x_j}^x \int_\eta^x d_{jki}(\xi,\eta,x,t)b_{jk}(\xi,\omega_j(\xi))
u_i^l(\eta,\omega_k(\eta;\xi,\omega_j(\xi))) d \xi d \eta\\
\displaystyle=\int_{x_j}^x \int_{\omega_j(\eta,x,t)}^{\omega_k(\eta,x,t)}
d_{jki}(\tilde x,\eta,x,t)\beta_{jk}(\tilde x,\omega_j(\tilde x))
\frac{(a_ka_j)(\tilde x,\omega_j(\tilde x))}{\partial_3\omega_k(\eta,\tilde x,\omega_j(\tilde x))}
u_i^l(\eta,\theta) d \theta d \eta,
\end{array}
\end{equation}
where the functions $\beta_{jk}\in C$ are fixed to satisfy (\ref{cass1}). 
Note that $\beta_{jk}(x,t)$ are not uniquely defined by (\ref{cass1}) for $(x,t)$ with 
$a_{j}(x,t)=a_{k}(x,t)$. Nevertheless, as it follows from (\ref{theta}), the right-hand side 
(and, hence, the left-hand side of (\ref{after_change}) do not depend on the choice of 
$\beta_{jk}$, since $d\theta=0$ if $a_{j}(\xi,\omega_j(\xi))=a_{k}(\xi,\omega_k(\xi))$.
Changing the order of integration in the right-hand side,
we rewrite the last as follows:
\begin{equation}\label{change_order}
\begin{array}{cc}
\displaystyle
 \int_{\omega_j(x_j)}^{\omega_k(x_j)}\int_{x_j}^{\tilde\omega_j(\theta)}
d_{jki}(\tilde x,\eta,x,t)\tilde b_{jk}(\tilde x,\omega_j(\tilde x))
\frac{(a_ka_j)(\tilde x,\omega_j(\tilde x))}{\partial_3\omega_k(\eta,\tilde x,\omega_j(\tilde x))}
u_i^l(\eta,\theta) d \eta d \theta
\\\displaystyle
+ \int^t_{\omega_k(x_j)}\int^{\tilde\omega_j(\theta)}_{\tilde\omega_k(\theta)}
d_{jki}(\tilde x,\eta,x,t)\tilde b_{jk}(\tilde x,\omega_j(\tilde x))
\frac{(a_ka_j)(\tilde x,\omega_j(\tilde x))}{\partial_3\omega_k(\eta,\tilde x,\omega_j(\tilde x))}
u_i^l(\eta,\theta) d \eta d \theta,
\end{array}
\end{equation}
where $\tau\in\mathbb{R}\mapsto \tilde\omega_s(\tau)=\tilde\omega_s(\tau,x,t)\in[0,1]$ 
denotes the inverse function to
 $\xi\in[0,1]\mapsto \omega_s(\xi)\in\mathbb{R}$.
Note that the interval of integration in $\theta$ in both of the integrals 
does not exceed $d$.
Now the $C\left(\overline\Omega_{2d}^{2d+\alpha}\right)^n$-norm
of the function (\ref{change_order}) can be estimated from above by
\begin{eqnarray}\label{conv1}
&2d\displaystyle
\max\limits_{x,\xi,\eta\in[0,1]}\max\limits_{t,\theta\in[0,2d+\alpha]}
\left|d_{jki}(\xi,\eta,x,t))\tilde b_{jk}(\xi,\theta)
\frac{\left(a_ka_j\right)
(\xi,\theta)}{\partial_3\omega_k(\eta,\xi,\theta)}\right|&
\nonumber\\
&\displaystyle\max\limits_{t\in[0,2d+\alpha]}\int_0^1|u_i^l(\eta,t)|\,d\eta
\displaystyle\le C\left\|u_i^l\right\|_{C([0,2d+\alpha],L^2(0,1))},&
\end{eqnarray}
where $C>0$ is a constant which depends on the coefficients of (\ref{eq:1}) but not on $u^l$.
The desired convergence of  $\left[B^2u^l\right](x,t)$ is thereby proved.

Now we treat the convergence of the sequence 
$$
\begin{array}{ll}
(RBu^l)_{j}(x,t)
=\displaystyle-c_{j}(x_j,x,t)\sum_{k=1}^n\sum_{s\neq k}\int_{0}^{1}\int_{x_k}^{\eta}r_{jk}(\eta,\omega_j(x_j))
d_{k}(\xi,\eta,\omega_j(x_j)) &  \\ [4mm]
\displaystyle \qquad\times b_{ks}(\xi,\omega_k(\xi,\eta,\omega_j(x_j)))u_{s}^l(\xi,\omega_k(\xi,\eta,\omega_j(x_j)))\,d\xi d\eta & 
\end{array}
$$
for an arbitrary fixed $j\le n$.
After changing the order of integration we get the equality
$$
\begin{array}{cc}
(RBu)_{j}(x,t)
=\displaystyle-c_{j}(x_j,x,t)\sum_{k=1}^n\sum_{s\neq k}\int_{0}^{1}\int_{\xi}^{1-x_k}r_{jk}(\eta,\omega_j(x_j))d_{k}(\xi,\eta,\omega_j(x_j)) &  \\ [2mm]
\displaystyle\times b_{ks}(\xi,\omega_k(\xi,\eta,\omega_j(x_j)))u_s^{l}(\xi,\omega_k(\xi,\eta,\omega_j(x_j)))\,d\eta d\xi.  & 
\end{array}
$$
Then we change the variable $\eta$ to $z=\omega_k(\xi,\eta,\omega_j(x_j))$. 
Since the inverse is given by $\eta=\tilde{\omega}_k(\omega_j(x_j),\xi,z)$,  we get
\begin{eqnarray}
\lefteqn{
(RBu)_{j}(x,t)
=-\displaystyle c_{j}(x_j,x,t)}\nonumber\\ &&\times\sum_{k=1}^n\sum_{s\neq k}
\displaystyle\int_{0}^{1}\int_{\omega_j(x_j)}^{\omega_k(\xi,1-x_k,\omega_j(x_j))}r_{jk}(\tilde{\omega}_k(\omega_j(x_j),\xi,z),\omega_j(x_j))
\label{rdoc}  
 \\ [2mm]
&&\times\displaystyle d_{k}(\xi,\tilde{\omega}_k(\omega_j(x_j),\xi,z),\omega_j(x_j))b_{ks}(\xi,z)
\displaystyle \partial_{3}\tilde{\omega}_k(\omega_j(x_j),\xi,z)u_s^{l}(\xi,z)\,dz d\xi,\nonumber
\end{eqnarray}
The functions $\omega_j(\xi,x,t)$ and the kernels of the integral operators in (\ref{rdoc}) are $C^1$-continuous. Due to the convergence 
(\ref{seq1}), there exists a limit of the right-hand side in 
$C\left(\overline\Omega_{2d}^{2d+\alpha}\right)$, as desired.

Now we consider the operator $G$. Changing the variable $\xi$ to $z=\omega_j(\xi,x,t)$ in (\ref{f3014}), we get
$$
%\begin{equation}
\begin{array}{ll}
(Gu^l)_j(x,t)= \\\displaystyle
-\sum_{k=1}^{n}\int_{\omega_j(x_j)}^{t}\int_{0}^{\tilde{\omega}_j(z)}d_{j}(\tilde{\omega}_j(z),x,t)g_{jk}(y,z)a_j(\tilde{\omega}_{j}(z),z)u_{k}^l(y,z)\, dydz.
\end{array}
$$
%\end{equation}
Similarly to the above, the functions $\omega_j(x_j), \tilde{\omega}_j(z), d_j(\tilde{\omega}_j(z),x,t)$, and $a_j(\tilde{\omega}_j(z),z)$ are continuous in $x$ and $t$. This entails the desired estimate of the type (\ref{conv1}).

We further proceed with the operator $R^2$. For $j\le n$, $k\le n$, and $T>0$, define operators $R_{jk}\in\mathcal{L}(C(\overline\Omega^T))$ by
$$
(R_{jk}w)(x,t)=c_j(x_j,x,t)\int_{0}^{1}r_{jk}(\eta,\omega_j(x_j))w(\eta,\omega_j(x_j))\,d\eta.
$$
Fix arbitrary $j\le n$, $k\le n$, and $i\le n$. We prove the needed estimate
for  the operator 
$R_{jk}R_{ki}$;  similar estimate for all other operators contributing into the $R^2$ will follow 
from the same argument. Given $T>0$, 
introduce operators $P_j,Q_{jk} : C(\overline\Omega^T)\to C(\overline\Omega^T)$ by
\begin{eqnarray}
(P_jw)(x,t)
&=&c_{j}(x_j,x,t)\int_{0}^{1}w(\eta,t)\,d\eta, \\
(Q_{jk}w)(x,t)
&=&r_{jk}(x,\omega_j(x_j))w(x,\omega_j(x_j)).
\end{eqnarray}
Then we have 
$$
R_{jk}=P_jQ_{jk}, \;\;\; R_{ki}=P_kQ_{ki}
$$
and, hence
$$
R_{jk}R_{ki}=P_jQ_{jk}P_kQ_{ki}.
$$
We aim at showing the estimate for $P_jQ_{jk}P_k$, as this and the boundedness of $Q_{ki}$ will entail the desired estimate for $R_{jk}R_{ki}$. The operator $P_jQ_{jk}P_k$ reads
\begin{equation}\label{ghng}
\begin{array}{lc}
(P_jQ_{jk}P_{k}w)(x,t) = c_j(x_j,x,t)  \\ [2mm]
\displaystyle \quad\times\displaystyle \int_{0}^{1} r_{jk}(\xi,\omega_j(x_j,\xi,t))c_k(x_k,\xi,\omega_j(x_j,\xi,t))\int_{0}^{1}w(\eta,\omega_k(x_k,\xi,t))\, d\eta d\xi.  &
\end{array}
\end{equation}
Changing the variable $\xi$ to $z=\omega_k(x_k,\xi,t)$, we get
\begin{equation}\label{ghng1}
\begin{array}{lr}
(P_jQ_{jk}P_{k}w)(x,t) =  c_j(x_j,x,t)  \\ [2mm]
\displaystyle\quad\times\int_{\omega_k(x_k,0,t)}^{\omega_k(x_k,1,t)}\hskip-5mm r_{jk}(\tilde{\omega}_k(t,x_k,z),z)c_k(x_k,\tilde{\omega}_k(t,x_k,z),z)\\ [2mm]\times \displaystyle \int_{0}^{1}\partial_{3}\tilde{\omega}_k(t,x_k,z)w(\eta,z)\, d\eta dz,  &
\end{array}
\end{equation}
where
\begin{equation} \label{2star}
\partial_3\tilde{\omega}_{k}(\tau,x,t)=
 a_k(x,t)\exp{\int_{\tau}^{t}\partial_1a_k(\tilde{\omega}_k(\rho,x,t),\rho)\, d\rho}.
\end{equation}
Similarly to the above, the needed estimate for $P_jQ_{jk}P_k$ now immediately follows from the regularity assumptions on the coefficients of the original problem.
We are therefore finished with the operator $R^2$.

Returning
 back to (\ref{conv_all}), it remains to treat the operator $BR$. By the definitions of $B$ and $R$,
\begin{equation}\label{dro}
\begin{array}{cc}
(BRu)_{j}(x,t)
=-\displaystyle \sum_{k\neq j}\sum_{l=1}^{n}\int_{0}^{1}\int_{x_j}^{x}d_{j}(\xi,x,t)b_{jk}(\xi,\omega_j(\xi))c_{k}(x_k,\xi,\omega_j(\xi))
&  \\
\displaystyle \times r_{kl}(\eta,\omega_k(x_k,\xi,\omega_j(\xi)))u_{l}(\eta,\omega_k(x_k,\xi,\omega_j(\xi)))\,d\xi d\eta, \;\;\; j\le n.    & 
\end{array}
\end{equation}
The integral operators in (\ref{dro}) are similar to those considered for $B^2$ and, therefore,
the proof follows along the same line as the proof for $B^2$. The proof of the theorem  is therefore complete.

\end{document}